\documentclass[12pt]{article}

\usepackage{amsmath}
\usepackage{upgreek}
\usepackage{enumitem} 
\usepackage{amsthm}
\usepackage{graphicx}
\usepackage{xcolor}
\usepackage{esdiff}
\usepackage{mathtools}
\usepackage{amssymb}
\usepackage{ marvosym }
\usepackage{amsmath,amssymb,amsfonts}
\usepackage{algorithmic}
\usepackage{graphicx}
\usepackage{textcomp}
\usepackage{xcolor}
\usepackage{float}
\usepackage[
backend=biber,
style=alphabetic,
sorting=ynt
]{biblatex}
\DeclareUnicodeCharacter{0308}{HERE!HERE!}
\bibliography{reference.bib}

\begin{document}

\title{Data Clustering and Visualization with Recursive Goemans-Williamson MaxCut Algorithm}
\author{An Ly, Raj Sawhney, Marina Chugunova}
\date{December 12th, 2023}
\maketitle

\begin{abstract}
In this article, we introduce a novel recursive modification to the classical Goemans-Williamson MaxCut algorithm, offering improved performance in vectorized data clustering tasks. Focusing on the clustering of medical publications, we employ recursive iterations in conjunction with a dimension relaxation method to significantly enhance density of clustering results. Furthermore, we propose a unique vectorization technique for articles, leveraging conditional probabilities for more effective clustering. Our methods provide advantages in both computational efficiency and clustering accuracy, substantiated through comprehensive experiments.
\end{abstract}

\section{Introduction}
Over the years, researchers have looked into the clustering of publications based on citations, textual similarities, and several other measures as mentioned in \cite{comparingrelatednessmeasures} or \cite{clusterbycitation}. Content-based clustering of articles and documents is a crucial task in the area of information retrieval and analysis, especially in the biomedical field. With over two million biomedical publications in existence, and more being published every year, several similarity based approaches have been attempted to accurately cluster these articles into groups as discussed in \cite{twomillionbiomedicalpublications} and \cite{contentbasedpublicationalgorithm}. 

MaxCut algorithms are wellknown and widely used to cluster vectorized datasets. It is worth noting that MaxCut optimization, as extensively discussed in  \cite{computation8030075}, \cite{GWA}, \cite{randomheuristics}, \cite{computationalapproach}, \cite{npcomplexity} is NP-complete. A detailed review of different MaxCut algorithms such as  semidefinite programming, a random strategy, genetic algorithms, combinatorial algorithms, and a divide-and-conquer strategy can be found in \cite{algorithmcomparison}.  Applications of  greedy randomized adaptive search procedures, variable neighborhood searches, and path-relinking intensification heuristics are studied in \cite{randomheuristics}. Some interesting new computational approaches for MaxCut approximation were analysed  in \cite{computationalapproach}. 

In this article, we explore recursive application of the Goemans-Williamson MaxCut Algorithm (GWA) combined with a higher-dimensional relaxation method. To the best of our knowledge it is not well known how the efficiency of clustering depends on the choice of the relaxation dimension in the original GWA. 

To ensure that medical publications with similar content are grouped together and are distinct from other clusters, our primary goal is to maximize the dissimilarity between these clusters using an unsupervised learning model. To do so, we introduce a weight matrix whose entries denote the similarity between medical publications by extracting keywords such as in \cite{unsupervisedkeyphraseextraction}. In our case, these weights are the distances between vectorized publications. The overall goal is a binary clustering; we want to split a given set of medical publications into two clusters, $A$ and $A^c$, such that the total dissimilarity between these two clusters is maximized. 

\section{A summary of GWA}
\label{GWASummary}
In general for a MaxCut algorithm, given a dataset we would like to split the index set (each index corresponds to an article in our case) to maximize the dissimilarity between the two resulting subsets.

In other words, given a matrix of similarity weights $W = \{w_{ij}\}$ where $i, j = 1, \cdots, n$, we want to split the index set $i = 1, \cdots, n$ into 2 sets $A$ and $A^c$ such that $\sum\limits_{i \in A, j \in A^c} w_{ij}$ is maximized. The expected value of total dissimilarity cut obtained by GWA will be at least 87.8\% of the optimal MaxCut value. The proof is detailed below for the purpose of showing that dimensions in the relaxation method can be chosen arbitrarily (the original proof can be found in \cite{GWA}).

First, let us introduce vector $y \in \{-1, 1\}^n$ such that $y_i = 1$ if $i \in A$ and $y_i = -1$ if $i \in A^c$. The original index clustering problem is now reduced to the following optimization problem: $\max \frac{1}{2} \sum\limits_{i < j} w_{ij}(1 - y_iy_j)$ subject to the constraint $y \in \{-1, 1\}^n$.

For the relaxation method that had been developed in \cite{2rankrelaxation}, we replace the scalars $y_i$'s with $n$-dimensional vectors $v_i \in S_n$ where $S_n$ is defined to be the $n$-dimensional unit sphere. In general, one can use dimensions that are higher than $n$ for the relaxation method, and we will explore efficiency of these higher dimensional relaxations numerically. We now replace the product of two scalars $y_iy_j$ with a dot product of two vectors $\langle v_i, v_j \rangle$.
The optimization problem then becomes: 
$\max \frac{1}{2} \sum\limits_{i<j} w_{ij}(1-\langle v_i,v_j\rangle)$ subject to the constraints $v_i \in S_n$ for all $i = 1, \cdots, n$. After some expansion of this expression, we can rewrite the above relaxed optimization problem as:
$\max \frac{1}{2} \sum\limits_{i<j} w_{ij} - \frac{1}{2} \sum\limits_{i<j} w_{ij} \langle v_i,v_j \rangle$. Note that in the above optimization problem, the first term is a fixed given scalar as all the weights are already predefined. 
Therefore, we can simplify the expression:
$\min \frac{1}{2}\sum\limits_{i<j} w_{ij}\langle v_i,v_j \rangle$. 

Now, let us introduce a matrix  $V = \{v_1 | v_2 | \cdots | v_n\}$. With some linear algebra, we can transform our optimization problem into the trace minimization problem:
$\min \frac{1}{4}\text{tr}(WV^TV)$ where $W$ is the original weight matrix and $V$ is defined as above. 
Let us define a new matrix $X = V^TV$. Note that all diagonal entries of positive semi-definite matrix $X$ will be equal to 1. 
Therefore, we can further rewrite this optimization problem as:
$\min \frac{1}{4}\text{tr}(WX)$ with the only constraint being that $v_i$ are vectors on the $n$-dimensional unit sphere for all $i$. 

Once we solve the optimization problem and find $V$, we will need to determine how our datapoints will get clustered and what will be the total dissimilarity value corresponding to this clustering. In order to cluster the datapoints into two groups, we will construct a separating hyperplane whose normal vector is $r$, a uniformly randomly chosen vector on our $n$-dimensional unit sphere $S_n$. 
To do so, we let $r$ be a vector that is randomly chosen from a multivariate normal distribution with the zero mean vector and the $n$-dimensional identity variance matrix, i.e. $r \sim N(0, I_n)$. 
After enough iterations and by normalizing random vectors to have unit lengths, we will generate a uniformly random vector on the $n$-dimensional plane. 

Now we have a separating hyperplane that splits our dataset of articles (recall each $v_i$ vector corresponds to an article) into two clusters. Let $C$ be the total dissimilarity value that corresponds to this cut. Let us define $\theta_{ij}$ to be the angle between vectors $v_i$ and $v_j$ for all pairs $i$ and $j$ in our index set.
Hence, we have $\cos{\theta_{ij}} = \langle v_i, v_j \rangle$ which will be used in later computations. Also note that, if $v_i$ and $v_j$ are on the same side of our hyperplane as determined by the dot product sign with the normal vector $r$, then the sum contribution term from this pair of indices will actually be equal to zero. 
\begin{equation*}
\begin{split}
    C = \frac{1}{2}\sum_{i <j} w_{ij}[1 - \langle v_i, v_j \rangle] \\
    \approx \frac{1}{2}\sum_{i < j} w_{ij} [ 1 - \text{sgn}(\langle v_i, r\rangle )\text{sgn}(\langle v_j,r\rangle )]    
    \end{split}
\end{equation*}
Using the properties of expected values, we will have:
\begin{align}
    E[C] & = E[\frac{1}{2}\sum_{i < j} w_{ij} [ 1 - \text{sgn}(\langle v_i, r\rangle )\text{sgn}(\langle v_j,r\rangle )]]\\
     & = \sum_{i < j} w_{ij} \text{Pr}(\text{sgn}(\langle v_i, r\rangle ) \neq \text{sgn}(\langle v_j,r\rangle ) \\
    & = \sum_{i < j} w_{ij} [2\text{Pr}(\langle v_i,r\rangle \geq 0, \langle v_j,r\rangle < 0)]
\end{align}
Note that this probability is proportional to the angle between vectors $v_i$ and $v_j$ as discussed in \cite{ComplexTheoryLecture} and \cite{SDPReview}. Specifically, this probability is equal to $\frac{\theta_{ij}}{2\pi}$. Continuing with this, we have:
\begin{align}
    E[C] & = \sum_{i < j} w_{ij} 2(\frac{\theta_{ij}}{2\pi})
    = \sum_{i < j} \frac{\theta_{ij}w_{ij}}{\pi}\\
    & = \sum_{i < j} \frac{\arccos{\langle v_i, v_j\rangle }w_{ij}}{\pi}
\end{align}
Recall that $w_{ij}$ will actually be counted as $0$ if the two vectors $v_i$ and $v_j$ are on the same side of the hyperplane determined by the sign of the dot product of $v_i$ and $v_j$ with the normal vector $r$. 

From here, we can determine the minimum value this expectation can take. Specifically, let us look at the expression $\theta_{ij}/\pi$. 
\begin{align*}
    \frac{\theta_{ij}}{\pi} & = \frac{\theta_{ij}}{\pi} \frac{2(1-\cos(\theta_{ij}))}{2(1-\cos(\theta_{ij}))}\\
    & = \frac{2\theta_{ij}}{\pi (1-\cos(\theta_{ij}))} \frac{1-\cos(\theta_{ij})}{2}
     \geq \alpha \frac{1}{2}(1-\langle v_i,v_j\rangle )
\end{align*} 
where we define $\alpha = \underset{0 < \theta \leq \pi}{\min} \frac{2\theta}{\pi(1-\cos(\theta))}$.

Therefore, the expected value of the total dissimilarity for this cut has a lower bound of:
$E[C] = \sum\limits_{i < j} \frac{w_{ij}\theta_{ij}}{\pi} \geq \alpha \frac{1}{2}\sum\limits_{i < j} w_{ij}(1-\langle v_i,v_j\rangle )$.
Note that the second part of this expression is exactly the solution of our relaxed optimization problem. 

It has been claimed and proven that our defined variable $\alpha > 0.878$ as discussed in \cite{GWAalpha} and in \cite{Feige1995ApproximatingTV}. 
This minimization is equivalent to minimizing the expression $\frac{\theta}{1-\cos(\theta)}$ under the same constraints.
This minimization occurs when $\theta = 2.331122\cdots$.

For a more formal proof as seen in \cite{GWA}, note that if we define $f(\theta) = 1-\cos(\theta)$, then $f(\theta)$ is concave for the interval $\frac{\pi}{2} < \theta < \pi$, i.e. we have that $f(\theta) \leq f(\theta_0) + (\theta - \theta_0)f'(\theta_0)$.
Replacing $f$ with it's definition, we have that:
\begin{align*}
    & 1 - \cos(\theta) \leq 1-\cos(\theta_0) + (\theta - \theta_0)\sin(\theta_0)\\
    & = \theta \sin(\theta_0) + (1-\cos(\theta_0) - \theta_0 \sin(\theta_0))
\end{align*}
where the last three terms create an overall negative term on our given interval.
When $\theta_0 = 2.331122\cdots$, then we have: $1 - \cos(\theta) < \theta \sin(\theta_0)$. 
Therefore, our function for the $\alpha$ in question follows the inequalities:
\begin{equation*}
    \frac{2\theta}{\pi(1-\cos(\theta))}  > \frac{2\theta}{\pi\theta\sin(\theta_0)}\\
     = \frac{2}{\pi\sin(\theta_0)} = 0.878\cdots 
\end{equation*}

Further improvements on this lower bound for the maxcut gained by any algorithm have been made using other maxcut algorithms such as the MAX $k$-Cut and Max Bisection algorithms with a lower bound between 0.800 and 0.990 depending on the $k$ used as studied in \cite{maxkcut}. Another study \cite{.699approx} has received a 0.699 approximation algorithm for Max-Bisection as well in comparison to the standard 0.50 approximation for the maxcut value.

Instance-specific lower bounds for the value of the cut gained from these algorithms were studied and compared across different algorithms such as randomized rounding, $k$-means, $k$-medoids, fuzzy $c$-means, and minimum spanning trees in \cite{computation8030075}, showcasing how clustering algorithms yield better cut values and improvements over the standardized randomized rounding algorithm. 

\section{Recursive Generalization of GWA}
During the optimization process we compute a set of $v_i$ vectors for each optimal matrix $X$ as discussed in Section~\ref{GWASummary}. If we were to use a uniformly random separating hyperplane method on the original dataset, the expected value of the cut would be just $50 \%$ of the optimal cut value. Using this uniformly random splitting method on the set of $v_i$ vectors results into an expected value that is higher than $87 \%$ of the optimal cut. That means $v_i$ vectors are clustered with better density and separation than the original dataset.   

In order to better understand this improvement of clustering due to mapping the original dataset to the set of $v_i$ vectors, we will analyse the recursive iterations of GWA. Specifically, we look to replace our initial dataset with the resulting $v_i$ vectors for further GWA iterations.

For this purpose, we use the $V$ matrix computed by the Cholesky Decomposition, and we also incorporate the principal component analysis of this $V$ matrix to visualize the results. 

Utilizing the Cholesky Decomposition, we construct the new dataset by treating each column of the $V$ matrix as a datapoint.
Due to the $n$-dimensionality of each $v_i$, good visualization is not possible without dimensionality reduction.

Hence, in order to visualize these vectors, we conduct PCA on $V$. In this article, we use 2-dimensional PCA as our dataset consists of the conditional probabilities of finding 'human' or 'side-effect' in the presence of amodiaquine as later discussed in Section~\ref{vectorization}. 

For our dataset, we would like to note that the effect of recursive GWA decreases with each iteration. Due to this, we only use up to an additional four iterations as, after that, no significant improvement is noticed. 

\subsection{Results of Recursive GWA}
\label{recursiveresults}
In this section, we now analyze the results of recursive GWA on two hand-generated datasets, the first dataset comprised of two separated cubes and the second dataset of interlocking clusters. Note that each of these datasets contains 100 points in total. 

\begin{figure}
\centering
        \includegraphics[width = 7.5cm]{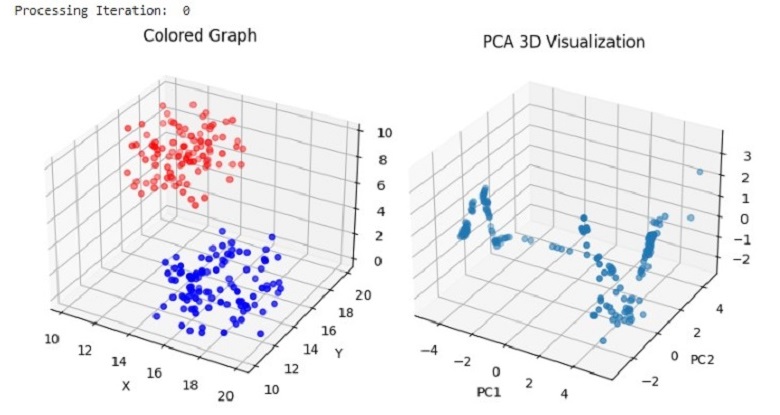}\\
        \includegraphics[width = 7.5cm]{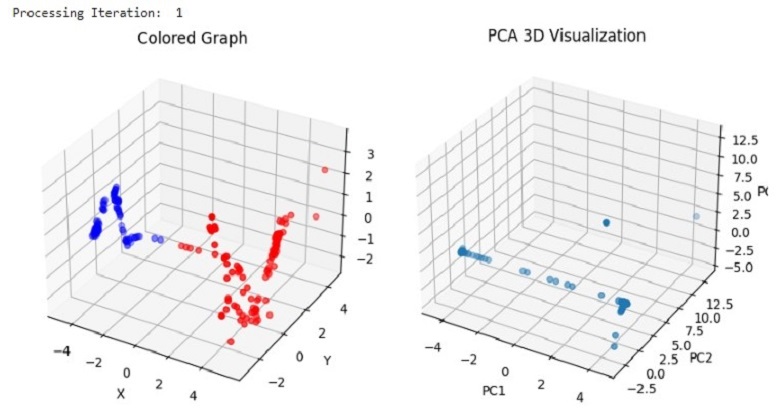}\\
        \includegraphics[width = 7.5cm]{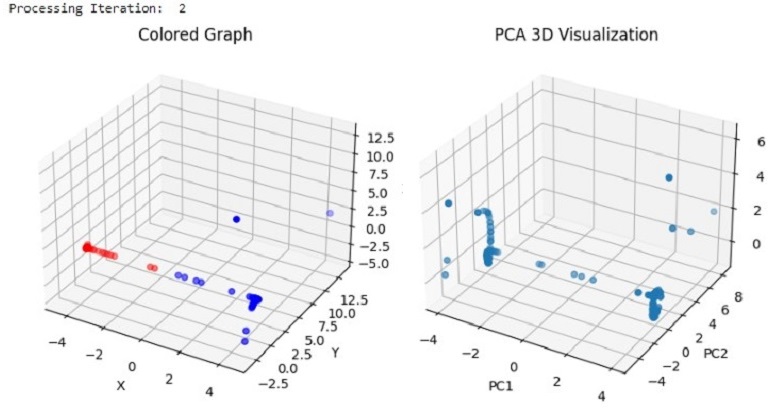}
    \caption{GWA Results and PCA Visualization over 3 Iterations: Two Cubes}
    \label{fig:iterativecube}
\end{figure}

Figure~\ref{fig:iterativecube} shows the results of our first dataset of two separated cubes. 
On initial iteration (top image), our model accurately separates the two cubes, and the $v_i$ vectors begin to cluster into two separate groups. 
Then, by using the 3D PCA results as input to the next iteration (middle image), we see tighter clustering of the updated $v_i$ vectors.
This trend continues in our final iteration (bottom image) as shown by the changing values on the axes. 

\begin{figure}
\centering
        \includegraphics[width = 7.5cm]{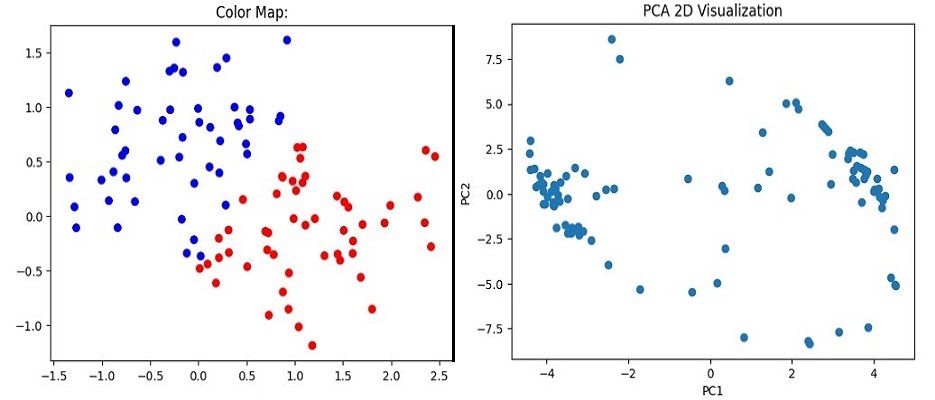}\\
        \includegraphics[width = 7.5cm]{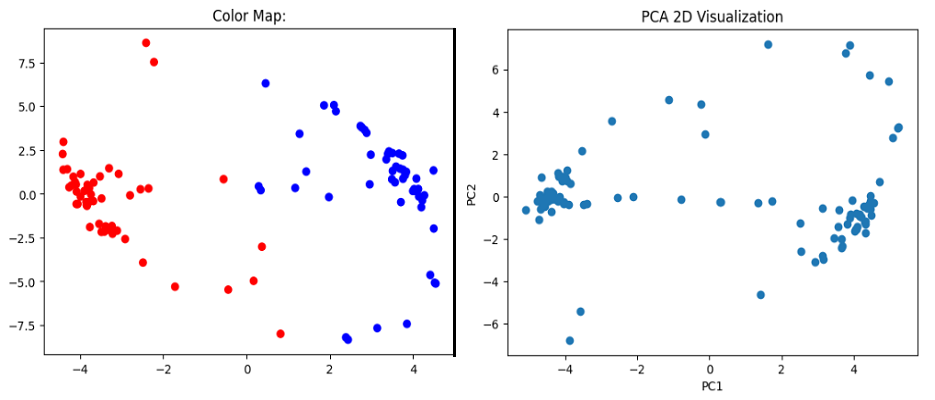}\\
        \includegraphics[width = 7.5cm]{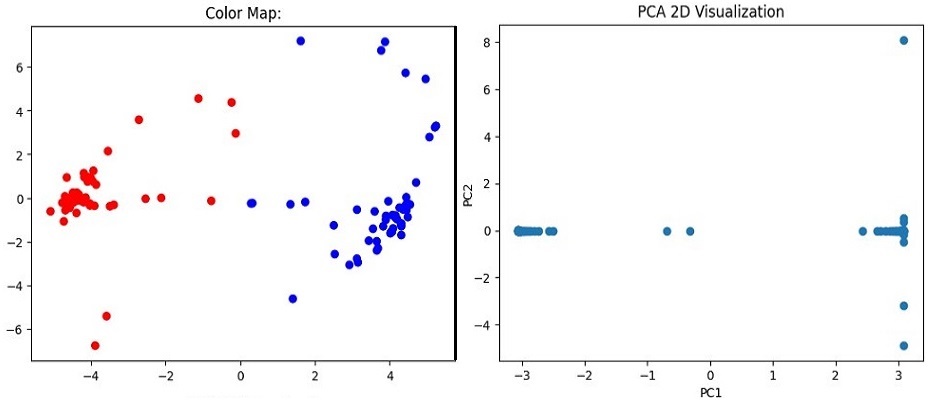}  
    \caption{GWA Results and PCA Visualization over 3 Iterations: Interlocking Data}
    \label{fig:iterativemoon}
\end{figure}

Figure~\ref{fig:iterativemoon} shows the results of using an interlocking dataset generated by the moonset command in Python. 
On initial iteration (top image), our model linearly separates the two interlocking groups, and the $v_i$ vectors begin clustering into two larger sets. 
By using the 2D PCA results as input for the next iteration (middle image), we see a slightly tighter clustering in the updated $v_i$ vectors.
In the final iteration of our data (bottom image), we see a significant clustering.

\section{Higher Dimensional Generalization of GWA}
In this section, we describe the mapping of the original dataset  into $m$-dimensional $v_i$ vectors where the dimension $m > n$ (recall that $n$ is the total number of points in the our original dataset).

We implement this higher-dimansional mapping by adding rows and columns of zeroes to the original weight matrix $W$ to match a new $m$ by $m$ dimension. 
Note that, because these values are zero, no new information is added to our dataset.

\subsection{Results of Higher Dimensional Generalization of GWA}
In this section, we now analyze the results of higher dimensional generalization on the dataset of interlocking clusters as first introduced in Section~\ref{recursiveresults}. 

\begin{figure}
\centering
        \includegraphics[width = 4cm]{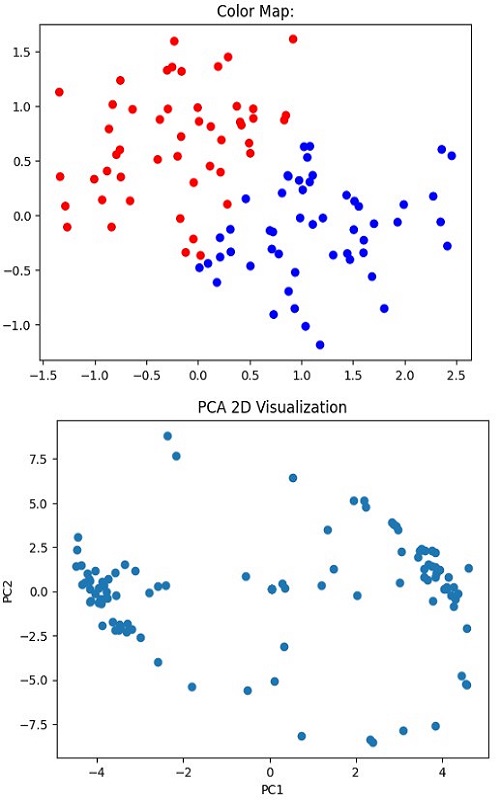}
        \includegraphics[width = 4cm]{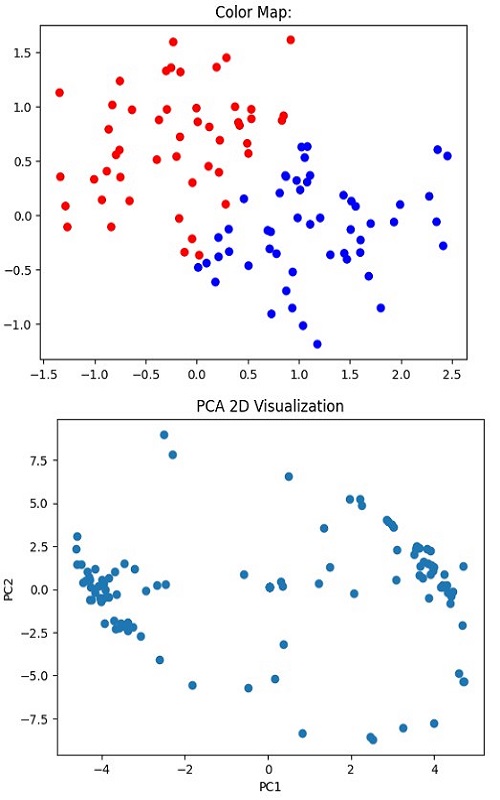}
    \caption{First Iteration with 104 Dimensions (Left) and 109 Dimensions (Right)}
    \label{fig:first104-109}
\end{figure}

Comparing Figure~\ref{fig:first104-109} with the results of Figure~\ref{fig:iterativemoon} (top image), no significant difference in clustering is seen. 
We now proceed to utilize both the higher dimensional generalization and the recursive generalization of GWA to analyze the effects on our model. 

\begin{figure}
\centering
        \includegraphics[width = 4cm]{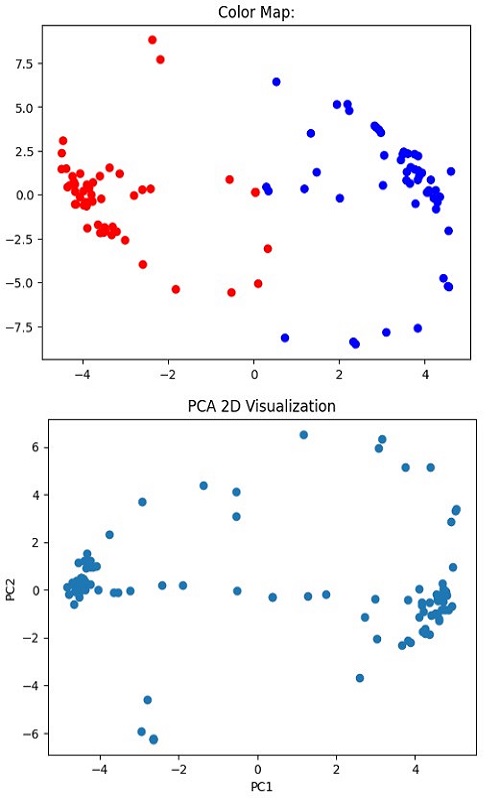}
        \includegraphics[width = 4cm]{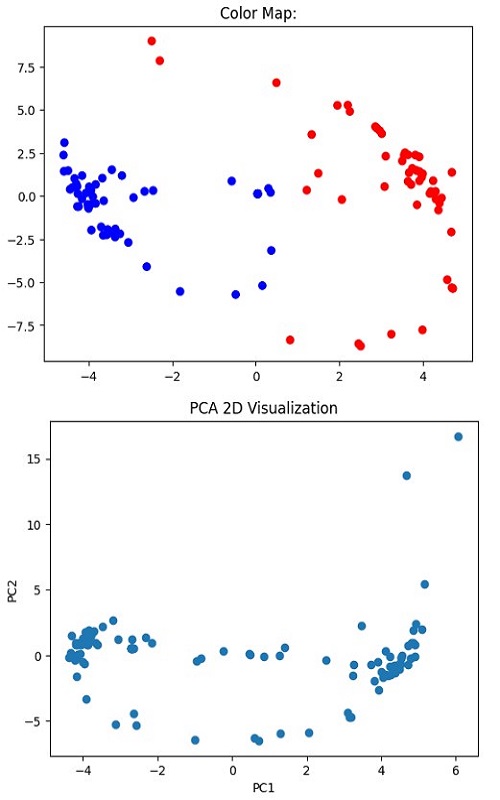}
    \caption{Second Iteration with 104 Dimensions (Left) and 109 Dimensions (Right)}
    \label{fig:second104-109}
\end{figure}

Using our 2D PCA visualizations from the initial iterations, we now compare Figure~\ref{fig:second104-109} and the middle image from Figure~\ref{fig:iterativemoon}. 
In comparing the results from all three cases, there is not a significant difference in clustering.
Excluding the outliers, the main clusters are tighter in 109 dimensions in comparison to 100 or 104 dimensions. 

\begin{figure}
\centering
        \includegraphics[width = 4cm]{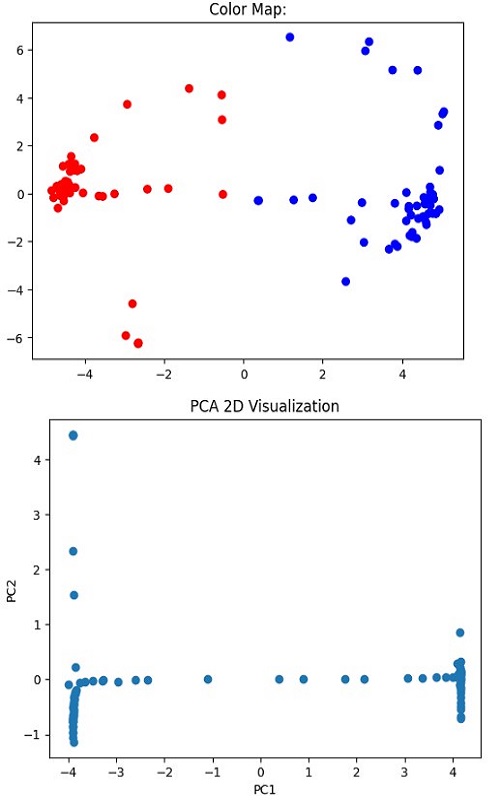}
        \includegraphics[width = 4cm]{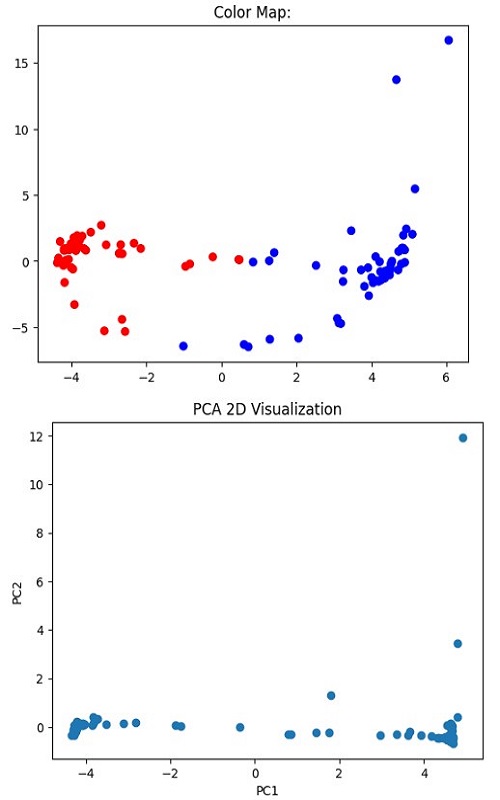}
    \caption{Third Iteration with 104 Dimensions (Left) and 109 Dimensions (Right)}
    \label{fig:third104-109}
\end{figure}

This trend continues in our final iteration as seen by comparing Figure~\ref{fig:third104-109} with the bottom images of Figure~\ref{fig:iterativemoon}. 
Our clusters in 100 dimensions are much tighter in comparison to 104 or 109 dimensions while still including outliers. 
While the results in 109 dimensions are worse than those in 100 dimensions, they are significantly better than the results in 104 dimensions.

\section{Vectorization of Articles} 
\label{vectorization}
We transform each article into a vector through the following procedure. We begin by identifying a list of target words such that for each given word, we compute the conditional probability of seeing the other target words within a set window size. As an initial approach, we define the following target list $L$: $L = [\text{amodiaquine}, \text{human}, \text{side-effect}]$.

Note that we first replace any medical side effect seen within the text, i.e. "headache", with the word "side-effect" before computing conditional probabilities. A similar procedure is used for human context words. The list of side effects is taken from \cite{sideeffectsebl}. These pre-processing steps allow for accurate probability vector calculations.

With this list, we set an arbitrary window of size $n$ centered around the target word. Within this window, we compute the conditional probabilities of seeing the other two target words. We repeat this procedure for every article and can therefore use these vectors to split articles of interest from irrelevant articles. We treat the window size $n$ as a hyper-parameter and, after testing, use $n = 10$ in the article "vectorization" process.

\section{Clustering Results of GWA}

\begin{figure}
\centering
        \includegraphics[width = 4.35cm]{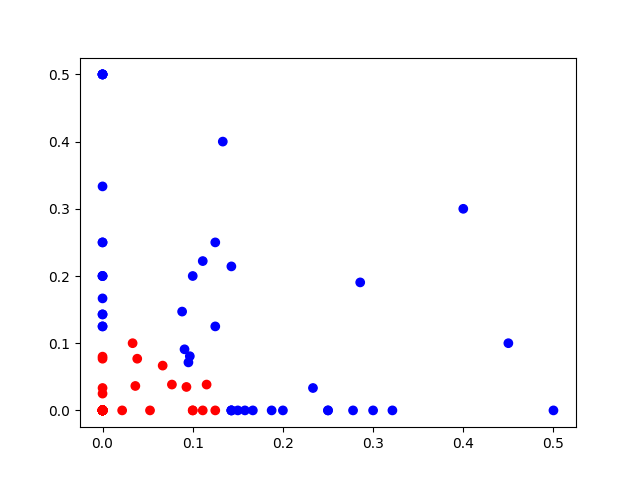}
        \includegraphics[width = 4.35cm]{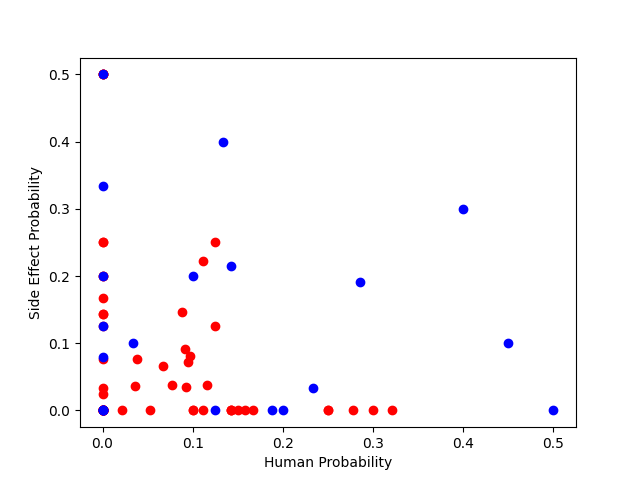}
    \caption{Comparison of Clustering Algorithms: GWA Clustering (Left) and True Labels as Humanly-Classified (Right) with a window size equal to 10 (5 on each side).}
    \label{fig:amodiaquine}
\end{figure}

Figure~\ref{fig:amodiaquine} shows the results of our initial iteration of GWA in comparison to the true labels. 
For the graph showcasing the datapoints and their true labels on the right, articles that discuss the side effects of amodiaquine in human patients are blue while those that do not are red. 
The horizontal axis denotes the probability that we find a 'human' context word in the presence of 'amodiaquine' while the vertical axis denotes the probability of finding a 'side-effect' in the presence of 'amodiaquine.' 
Note that some of our desired articles are at $(0,0)$ meaning that there is no mention of a side effect or a human patient within the presence of 'amodiaquine'. 
This is most likely due to the window size being too small.

The left image in Figure~\ref{fig:amodiaquine} shows the results of GWA clustering on our medical paper dataset. 
Note that we get a separation line that is almost completely linear between the two clusters. 
This algorithm manages to cluster almost all of the true articles into one grouping. 
The existing outliers are likely due to the placement of these $(0,0)$ articles. 

\section{Conclusions}
In conclusion, our article presents modified Goemans-Williamson MaxCut Algorithms  and studies their applications to clustering vectorized datasets. Namely, we demonstrate the effectiveness of recursive iterations and higher-dimensional generalizations of the GWA in achieving dissimilarity-based clustering. We think these methods combined with dimensionality reductions have the potential to further enhance clustering results. In addition, the introduction of the vectorization method based on conditional probabilities provides an additional tool for unsupervised document classifications.

While GWA shows promise in accurately clustering articles, there are some challenges that will need to be researched and refined. Future development of techniques to handle outliers and to fine-tune our parameters will contribute to a more precise and robust method.

\printbibliography
\end{document}